\documentclass[12pt,a4paper,leqno]{article}

\usepackage{latexsym}
\usepackage{amssymb,exscale}
\usepackage[centertags]{amsmath}
\usepackage{amsthm}
\usepackage[all]{xy}

\numberwithin{equation}{section}

\swapnumbers
\newtheorem{theorem}[equation]{Theorem}
\newtheorem*{maintheorem}{Main Theorem}
\newtheorem*{mainlemmanonumber}{Main Lemma}
\newtheorem{mainlemma}[equation]{Main Lemma}
\newtheorem{lemma}[equation]{Lemma}
\theoremstyle{definition}

\renewcommand{\phi}{\varphi}

\renewcommand{\(}{\bigl(}
\renewcommand{\)}{\bigr)\vphantom{)}}

\newcommand{\Lip}{\operatorname{Lip}}
\newcommand{\N}{\operatorname{N}}
\newcommand{\const}{\operatorname{const}}
\newcommand{\Var}{\operatorname{Var}}
\newcommand{\dist}{\operatorname{dist}}

\newcommand{\Const}{\operatorname{Const}}

\newcommand{\pd}{\partial}

\newcommand{\K}{\mathrm K}
\newcommand{\Om}{\Omega}
\newcommand{\om}{\omega}

\newcommand{\eps}{\varepsilon}
\newcommand{\ga}{\gamma}

\newcommand{\la}{\lambda}
\newcommand{\de}{\delta}
\newcommand{\De}{\Delta}
\newcommand{\F}{\mathcal F}
\newcommand{\cN}{\mathcal N}

\newcommand{\Ex}{\mathbb E}
\renewcommand{\Pr}[1]{\,\mathbb P\,\(\,#1\,\)\,}
\newcommand{\R}{\mathbb R}
\newcommand{\C}{\mathbb C}
\newcommand{\Z}{\mathbb Z}

\newcommand{\ti}{\tilde}
\newcommand{\sif}{$\sigma$-field}

\def\One{{1\hskip-2.5pt{\rm l}}}

\def\emailwww#1#2{\par\qquad {\tt #1}\par\qquad {\tt #2}\medskip}

\begin{document}

\title{Random complex zeroes, II. \\
 Perturbed lattice}

\author{Mikhail Sodin\thanks{Supported by the Israel
Science Foundation of the Israel Academy
of Sciences and Humanities}
\ and Boris Tsirelson}

\date{}

\maketitle

\begin{abstract}
We show that the flat chaotic analytic zero points
(i.e. zeroes of a random entire function $ \psi(z) = \sum_{k=0}^\infty
\zeta_k \frac{z^k}{\sqrt{k!}} $
where $\zeta_0, \zeta_1, \dots $ are independent standard
complex-valued Gaussian variables)
can be regarded as a random perturbation of a lattice in the plane. The
distribution of the distances between the zeroes and the corresponding lattice
points is shift-invariant and has a Gaussian-type decay of the tails.
\end{abstract}

\section*{Introduction}

We consider the (random) set $ S $ of zeroes of a random entire
function $ \psi \colon \C \to \C $,
\begin{equation}\label{0.0}
\psi(z) = \sum_{k=0}^\infty \zeta_k \frac{z^k}{\sqrt{k!}} \, ,
\end{equation}
where $\zeta_0, \zeta_1, \dots $ are independent standard
complex-valued Gaussian random variables; that is, the
distribution $\mathcal N_{\C}(0,1)$ of each $\zeta_k$ has the
density $\pi^{-1} \exp(-|w|^2)$ with respect to the Lebesgue
measure $m$ on $\C$. Well-known as the flat CAZP (`chaotic analytic
zero points'), this model is distinguished by invariance of
the distribution of zero points with respect to the motions of the
complex plane, see \cite{SoTsi} for details and references.

\subsubsection*{Toy models}

It is instructive to compare the flat CAZP with simpler (`toy') models
of random point processes in the plane, especially, random
perturbations of a lattice. The first toy model: each point of the
lattice $ \sqrt{\frac\pi2} \Z^2 = \{ \sqrt{\frac\pi2} (k+li) \colon k,l \in
\Z \} $ is deleted at random, independently of others, with
probability $ 1/2 $; the remaining points are a random set $ S_1
$. For smooth functionals (linear statistics)
\[
Z_{L,h} (S) = \sum_{z\in S} h \bigg( \frac z {\sqrt L} \bigg) \, ,
\]
where $ h \colon \C \to \R $ is a compactly supported smooth function, mean
values are similar for $ L \to \infty $,
\begin{align*}
\Ex Z_{L,h} (S) &\sim L \frac1\pi \int h \, dm \, , \\
\Ex Z_{L,h} (S_1) &\sim L \frac1\pi \int h \, dm \,
\end{align*}
(here and below, $m$ always stands for the Lebesgue measure), but
fluctuations of $ S_1 $ are much stronger:
\begin{align*}
\Var Z_{L,h} (S) &\sim \frac{\const}L \| \De h \|^2 \, , \\
\Var Z_{L,h} (S_1) &\sim \const \cdot L \| h \|^2 \, ,
\end{align*}
see \cite{SoTsi}, the end of the introduction.

The second toy model: points of the lattice $ \sqrt\pi \Z^2 $ move
independently, giving
\[
S_2 = \{ \sqrt{\pi} (k+li) + \eta_{k,l} \colon\, k,l \in \Z \} \, ,
\]
where $ \eta_{k,l} $ are independent standard complex Gaussian random
variables. We have \cite{SoTsi}
\begin{align*}
\Ex Z_{L,h} (S_2) &\sim L \frac1\pi \int h \, dm \, , \\
\Var Z_{L,h} (S_2) &\sim \const \cdot \| \nabla h \|^2 \, ;
\end{align*}
the latter is closer (than $ L \| h \|^2 $) to $ L^{-1} \| \De h \|^2
$, but still dissimilar.

Asymptotic similarity to $ S $ can be reached (see the third toy model
in \cite{SoTsi}) by inventing special correlation between
perturbations $ \eta_{k,l} $.

\subsubsection*{Main result}

Discarding toy models and asymptotic properties, we come to the
idea of CAZP as a perturbed lattice,\footnote{%
 Area of its cells must equal $\pi$. Any lattice with this cell area may be
 used.}
\[
S = \{ \sqrt{\pi} (k+li) + \xi_{k,l} \colon\, k,l \in \Z \}
\]
for some (dependent) complex-valued random variables $ \xi_{k,l} $. Of
course, it can happen that all points of $ S $ are far from the
origin, in which case $ |\xi_{0,0}| $ is necessarily large. However,
that is an event of small probability. We may hope for fast decay of
the probability $ \Pr{ |\xi_{k,l}| \ge r } $ for large $ r $,
uniformly in $ k,l $. The uniformity becomes trivial if random
variables $ \xi_{k,l} $ are identically distributed. Taking into
account invariance of CAZP under shifts of $ \C $ we may hope for
invariance of $ (\xi_{k,l}) $ under lattice shifts. The hopes come
true, which is our main result, formulated below. Random variables are
treated as measurable functions on the space $ \Om $ of
two-dimensional arrays $ \xi \colon \Z^2 \to \C $ of complex numbers.

\begin{maintheorem}
There exists a probability measure $ P $ on
\textup{(}the Borel \sif\ of\textup{)} the space $ \Om = \C^{\Z^2}
$, invariant under shifts of $ \Z^2 $ and such that

\noindent\textup{(a)} the random set $ \{ \sqrt{\pi} (k+li) +
\xi_{k,l} \colon\, k,l \in \Z \} $ is distributed like the flat CAZP;

\noindent\textup{(b)} $ \Ex \exp \( \eps |\xi_{0,0}|^2 \) < \infty $
for some $ \eps > 0 $.
\end{maintheorem}

Item (a) needs some comments. A `random set' is a measurable map
from $ \Om $ to a space of sets. We need only locally finite
subsets of $ \C $. The Borel \sif\ on that space is generated by
functions $ S \mapsto \sum_{z\in S} h(z) $, where $ h $ runs over
compactly supported continuous (or just Borel) functions $ \C \to
\R $. Alternatively, we may represent each set $ S $ by its
counting measure, which is the same, since the random entire
function $ \psi $  has only simple zeroes (almost surely). Item
(a) means that the two maps
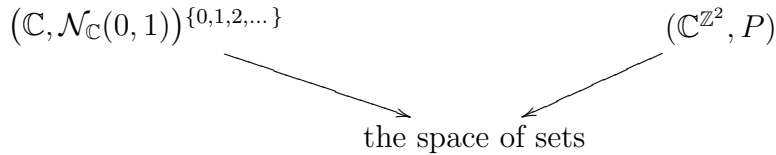
\[
\begin{gathered}
\xymatrix{
 \( \C, \mathcal N_\C(0,1) \)^{\{0,1,2,\dots\}} \ar[dr] & &
  (\C^{\Z^2},P) \ar[dl] \\
 & \text{the space of sets}
}
\end{gathered}
\]
induce the same measure on the space of sets. Here the first map sends
a sequence of coefficients $ \zeta_0, \zeta_1, \dots $ into the set of
zeroes of $ \psi(z) = \sum \zeta_k z^k / \sqrt{k!} $, while the second
map sends an array $ (\xi_{k,l})_{k,l} $ into the set $ \{ \sqrt{\pi}
(k+li) + \xi_{k,l} \colon\, k,l \in \Z \} $. The latter set is locally
finite (almost surely), which is a part of item (a).

The second map on the diagram intertwines natural (measure preserving)
actions of lattice shifts on $ \Om $ and the space of sets. For the
first map, the situation is more complicated; only a \emph{projective}
action of shifts is naturally defined on the space of entire functions
(or their coefficients), see \cite{SoTsi}.

Our construction of the matching between the flat CAZP and the
lattice points $ \sqrt{\pi}\Z^2 $ is not explicit. For this
reason, the main theorem gives no information about correlations
between $ \xi_{k,l} $. On large distances, the correlation
function of the underlying Gaussian process decays rapidly (see
for example \cite[Sect.~3.2]{SoTsi}). Probably, $ \xi_{k,l} $ can
be chosen as to be nearly independent on large distances. It is
also compatible with the result of \cite{SoTsi3}. On small
distances we expect a negative correlation, for two reasons: the
well-known repulsion of close zeroes \cite{Ha96, BSZ00a}, and the
center-of-mass conservation discussed in
\cite[Intoduction]{SoTsi}.

\subsubsection*{Explicit matching?}

It could be very useful to find an explicit matching between the
flat CAZP and the lattice points, or equivalently, a
transportation of the Lebesgue measure $ \frac1{\pi}m $ to the
counting measure $ n_\psi $ of CAZP. By `transportation' we mean a
map $ T\colon \C\to\C $ such that $ \pi n_\psi = T_* m $. Of
course, we are interested in stationary random transportations $ T $ with fast
decay of the tails of the distribution of $ Tz-z $. Here, we suggest a
natural and explicit construction of the transportation which, in
our opinion, deserves a better look\footnote{%
 It is inspired by the celebrated Moser homotopic construction
 \cite[Sect.~4]{Moser} of the diffeomorphism that transports one volume
 measure to another.}

Consider the gradient field of the stationary random potential
\cite[Introduction]{SoTsi}
\begin{equation}\label{potential}
\phi (z) = \frac12 \log|\psi(z)| - \frac14 |z|^2
\end{equation}
(additional factor $ \frac12 $ on the RHS will be convenient
later); the distributional Laplacian of $ \phi $ equals $ \pi
n_\psi - m $. The only local minima of the function $ \phi $ are
the points where it equals $ -\infty $ (since $ \phi $ is
superharmonic everywhere except of these points). We say that the
point $ w $ belongs to the \emph{basin} of a zero
point $ z\in\psi^{-1}(0) $ if $ \nabla\phi(w)\ne 0$, and the
gradient trajectory passing through $ w $ terminates at $ z $.
In other words, if we put a ball at the point $ (w,
\phi(w)) $ on the graph of the function $ \phi $, then under
the gravitation force (and without inertia) the ball will fall through at the
point $ (z, -\infty) $.

We expect that almost surely one obtains a cellulation of the
plane on \emph{finite cells,} each of them being the basin of some
random zero point; i.e.\ with probability one, nothing escapes to
infinity, and nothing arrives from infinity. On the boundary of
each cell, the gradient $ \nabla\phi $ has zero normal component.
Since each cell contains exactly one random zero point, by Green's
theorem applied to the function $ \phi $, the area of each cell
must equal $ \pi $. Define $ T $ as the map that sends each basin
into the corresponding zero point, then $ T $ transports the
measure $ \frac1{\pi}m $ to the measure $ n_\psi $. It would be
interesting to obtain a good estimate for the diameters of the
basins.

\subsubsection*{Other point processes}

It is instructive to think about possible counterparts of the main
theorem for simpler random sets such as the first toy model
(Bernoulli process) or the Poisson point process. Here, $ \Ex
|\xi_{0,0}| $ must be infinite, since for a finite fragment of
size $ n \times n $ the transportation cost between the random set
and the lattice, divided by the number of points, typically
exceeds $ \const \cdot \sqrt{ \ln n } $ \cite{AKT, T}. The large
gap between $ \Ex |\xi_{0,0}| $ and $ \Ex \exp \( \eps
|\xi_{0,0}|^2 \) $ shows that the random zeroes are distributed
much more evenly than independent random points (see \cite{PV}, Fig.~1
and comments to it). Note also that
(i) stationary matchings between random and deterministic sets are
closely related with so-called extra head schemes \cite{HP}; (ii)
existence of a stationary matching between the Poisson point
process on $ \R^2 $ and the lattice follows from existence of a
`stable marriage of Poisson and Lebesgue' announced recently
\cite{HP, HHP}.

\subsubsection*{The reader's guide, I (informal)}

The proof of the main theorem is based on the formula
\begin{equation}\label{0.3}
2\pi \, dn_\psi = \De \log |\psi| \, dm \, ,
\end{equation}
where $ n_\psi $ is the counting measure on the set of zeroes of $
\psi $. The desired array $ (\xi_{k,l})_{k,l} $ may be thought of
as a bijective correspondence (`marriage') between zeroes of $
\psi $ and lattice points. The distance $ |\xi_{k,l}| $ between
corresponding points (`fianc\'e' and `bride') must be controlled
as to ensure item (b) of the theorem. It is instructive to try
first a simpler condition, say, $ |\xi_{k,l}| \le 100 $ for all $
k,l $. (In fact, it is too much for a typical $ \psi $, but let us
try to prove it anyway.) If such $ \xi_{k,l} $ exist, then clearly
\begin{equation}\label{0.4}
n_\psi (U) \le n(U_{+r}) \quad \text{and} \quad n(U)
\le n_\psi (U_{+r})
\end{equation}
for every $ U \subset \C $; here $ U_{+r} $ stands for the $ r
$-neighborhood of $ U $, $ r=100 $, and $ n $ is the counting measure
on the
lattice. (No measurability is required of $ U $, since the measures $
n, n_\psi $ are discrete. However, it does not harm to assume $ U $ to
be a bounded domain with a smooth boundary.) In fact, \eqref{0.4} is
necessary \emph{and sufficient,} which is basically the well-known
`marriage lemma'. Using \eqref{0.4} as a sufficient condition, we may
replace $ n $ by $ \frac1\pi m $ at the expense of some change of the
constant $ r $.

Taking into account that $ \( \Ex |\psi(z)|^2 \)^{1/2} = \exp \(
\frac12 |z|^2 \) $ one could expect na\"{\i}vely
that the `potential' \eqref{potential}
is bounded on $ \C $. This can be used to show that $ \int_U \De\phi \, dm \le
m ( U_{+r} \setminus U ) $ and $ -\int_{U_{+r}} \De\phi \, dm \le m ( U_{+r}
\setminus U ) $, which gives \eqref{0.4} since $ \De \phi = \pi n_\psi - m $.
Singularity of the potential at zeroes is not an obstacle,
since we can replace $ n_\psi $ by its convolution with a compactly
supported smooth measure.

The argument sketched above does not work since the smoothed
random potential is unbounded (for almost all $ \psi $). However,
it can be mended. Rare fluctuations appear somewhere on the
infinite plane $ \C $, probably far from the origin. In order to
get item (b) of the theorem we need some locality; $ |\xi_{0,0}| $
should not be large, whenever the potential is not large in an
appropriate neighborhood of the origin. Restrictions on $
|\xi_{k,l}| $ should be adaptive, they should be relaxed around
large values of the potential. This idea is formalized by
introducing on $ \C $ a metric $ \rho $ that depends on $ \psi $,
and considering $ \rho $-neighborhoods $ U_{+r} $. Such metrics $
\rho $ are not shift-invariant; rather, the probability
distribution on the space of these metrics is shift-invariant.
Finally, shift invariance of $ P $ is ensured.

\subsubsection*{The reader's guide, II (more formal)}

To build the metric $ \rho $, we use an idea borrowed from
\cite[Section~1.4]{Hormander}. A $ \Lip(1) $-function $ R\colon
\R^d \to (0,\infty) $, $ |R(x)-R(y)|\le |x-y| $, gives rise to the
metric
\[
\rho (x, y) = \inf_\gamma  \int_\gamma \frac{|dz|}{R(z)}
\]
where the infimum is taken over all piece-wise $ C^1 $-curves $
\gamma $ in $ \R^d $ connecting the points $ x $ and $ y $. We
shall call $ \rho $ a {\em special metric} on $ \R^d $,
corresponding to $ R $ (in which case $ R $ is always assumed to
be a positive $ \Lip(1) $ function). In Sect.~1, we list the
properties of the metric $ \rho $ needed in the rest of the paper;
in Sect.~2, we construct a Whitney-type partition of unity
subordinated to the metric $ \rho $.

This partition of unity is needed for the following potential
theory lemma which lies in the heart of our argument:

\begin{mainlemmanonumber} There exists $ \const $ such that if a $ C^2
$-function $ u\colon \R^d \to \R $ and a special metric $ \rho $
\textup{(}corresponding to $ R $\textup{)} satisfy
$ |u(x)| \le \const \cdot \, R^2 (x) $ and $ \De u (x) \ge -1 $
for all $ x $, then
\[
\int_U \De u \, dm \le m ( U_{+4} \setminus U ) \, , \qquad -
\int_{U_{+4}} \De u \, dm \le m ( U_{+4} \setminus U )
\]
for every compact set $ U \subset \R^d $.
\end{mainlemmanonumber}
\noindent Here $ U_{+4} = \{ y \colon \exists x \in U \; \rho(x,y) \le
4 \} $. For the proof see Sect.~3.

\medskip
To apply this lemma, we smooth the random potential $ \phi $ (see
\eqref{potential})
setting $ u=\phi*\chi $, where $ \chi $ is a compactly supported
smooth function $ \chi \colon \C \to [0,\infty) $ such that $ \int \chi
\, dm = 1 $ and $ \chi(-z) = \chi(z) $ for all $ z $ (the choice
of $ \chi $ influences only constants), and define the $ \Lip(1)
$-function $ R $ as follows (see \eqref{4.4}):
\[
R(z) = \max_w \( \sqrt{ \Const \cdot \, ( 1 + (|\phi|*\chi)(w) ) } -
|w-z| \) \, .
\]
Then the main lemma combined with the marriage lemma and the general
inequality $ |x-y|\le 2^{3\rho (x,y)} R(x) $ (see Lemma \ref{1.5}) give us a
matching between CAZP and the lattice points $ \sqrt{\pi}\Z^2 $
such that for every matched pair $ z\in\psi^{-1}(0) $ and $
\sqrt{\pi}(k+il) $
\[
|z - \sqrt{\pi}(k+il)| \le R\left( \sqrt{\pi} (k+il)\right)
\]
(see Theorem~\ref{4.1}).

It will be shown (see Lemma \ref{5.1}) that $ \Ex \exp \( c R^2 (x) \) \le C
$. This will give us the subgaussian decay
of the tails, thus proving Item (b) of the main theorem.

It worth mentioning that the Gaussian nature of the random
function $ \psi $ is used only once, in Sect.~5, when checking the
inequality $ \Ex \exp ( \const \cdot |\phi(z)| ) \le \Const $
(uniformly in $ z $). For every random entire function satisfying
the inequality, its zeroes are a `perturbed lattice' satisfying $
\Ex \exp ( \const \cdot |\xi_{k,l}|^2 ) < \Const $ (uniformly in $
k,l $). On the other hand, shift-invariance of $ (\xi_{k,l})_{k,l}
$ is achieved via shift-invariance of zeroes of $ \psi $; the
latter was verified using the Gaussian distribution.

\subsubsection*{Convention}

Most of the steps in the proof of the main theorem do not use any
special properties of the complex plane and will be done in the
Euclidean space $ \R^d $. Throughout, `$ \Const $' and `$ \const
$' mean positive constants (sufficiently large and sufficiently
small, respectively) depending on the dimension $ d $ only, the
values of these constants can be changed at each occurrence. By $
B(x; r) $ we always denote the closed ball $ \{y\colon |x-y|\le r \} $.

\subsubsection*{Acknowledgment} We thank Michael Krivelevich for
providing us with the reference \cite{Aharoni} and its discussion,
F\"edor Nazarov for bringing to our attention the relevance of the
classical Whitney construction, and Yuval Peres for telling us
about results of the papers \cite{AKT, T, HHP, HP}.

\section{A class of metrics in $\R^d$}

Suppose $ \rho $ is the special metric corresponding to a $\Lip(1)
$-function $ R\colon \R^d \to (0,\infty) $; i.e.
\[
\rho (x, y) = \inf_\gamma  \int_\gamma \frac{|dz|}{R(z)}
\]
where the infimum is taken over all piece-wise $ C^1 $-curves $
\gamma $ in $ \R^d $ connecting the points $ x $ and $ y $. We
prove several simple facts about the special metric $ \rho $.

\begin{lemma}\label{1.2} For all $ x,y \in \R^d $,
\[
|y-x| \le \frac12 R(x) \quad \text{implies} \quad \frac12 R(x) \le
R(y) \le \frac32 R(x).
\]
\end{lemma}

\begin{proof} $ R(y) \le R(x) + |x-y| \le \frac32 R(x) $, and
$ R(x) \le R(y) + |x-y| \le R(y) + \frac12 R(x) $ which gives $
\frac12 R(x) \le R(y) $.
\end{proof}

Given a piece-wise $ C^1 $-curve $ \ga $ which starts at $ x $
and terminates at $ y $, we define the \emph{index} $ N(\ga)\in\N $
as the length $ N $ of the chain  of points $ x_0=x $, $ x_1 $, ..., $
x_N=y $ on $ \ga $ constructed one after another as follows. Having
$ x_j $, we consider the rest $ [x_j,y]_\ga $ of the curve $ \ga $
(the part of $ \ga $ which starts at $ x_j $ and terminates at $ y
$). If $ [x_j,y]_\ga \subset B \( x_j; \frac12 R(x_j) \) $ then the
process stops at $ x_{j+1} = y $, $ N=j+1 $. Otherwise $ x_{j+1} $ is
the first point on $ [x_j,y]_\ga $ lying on the sphere $ \pd B \( x_j,
\frac12 R(x_j) \) $, and the process is continued.

\begin{lemma}\label{1.3} $ N(\ga) \le 3 \int_\ga \frac{|dz|}{R(z)} +
1$.
\end{lemma}
\begin{proof} Denote by $ [x_j,x_{j+1}]_\ga $ the part of $ \ga $
between $ x_j $ and $ x_{j+1} $. If $ z\in [x_j,x_{j+1}]_\ga $, then $
R(z) \le \frac32 R(x_j) $ by Lemma \ref{1.2}, therefore, denoting $ N
= N(\ga) $ and assuming $ N > 1 $ (otherwise there is nothing to
prove),
\begin{multline*}
\int_\ga \frac{|dz|}{R(z)} = \sum_{j=0}^{N-1} \int_{[x_j,x_{j+1}]_\ga}
 \frac{|dz|}{R(z)} \ge \sum_{j=0}^{N-2} \int \dots \\
\ge \sum_{j=0}^{N-2} \frac{2}{3R(x_j)} \cdot |x_j-x_{j+1}| =
 \sum_{j=0}^{N-2} \frac{2}{3R(x_j)} \cdot \frac12 R(x_j) =
 \frac{N-1}{3}\, .
\end{multline*}
\end{proof}

\begin{lemma}\label{1.4} Let $ \rho $ be a special metric on $
\R^d $, corresponding to $ R $. Then
\[
|x-y|\le \frac12 R(x) \quad \text{implies} \quad  \rho (x,y) \le 1 \,
 .
\]
\end{lemma}

\begin{proof}
Let $ [x,y]\subset \R^d $ be the straight segment with end-points at $
x $ and $ y $. Then
\[
\rho (x,y) \le \int_{[x,y]} \frac{|dz|}{R(z)}
\stackrel{\ref{1.2}}\le \frac{|x-y|}{\frac12 R(x)} \le 1 \, .
\]
\end{proof}

\begin{lemma}\label{1.5} Let $ \rho $ be a special metric on $
\R^d $, corresponding to $ R $. Then
\[
|x-y|\le 2^{3\rho (x,y)} R(x) \, .
\]
\end{lemma}

\begin{proof} Given $ \eps>0 $,  choose a curve $ \ga $
connecting the points $x $ and $ y $,  such that $ \int_\ga
\frac{|dz|}{R(z)} < \rho (x,y) + \eps $. Let $x_0=x$, $x_1$,
..., $x_N=y$ be a partition of $ \ga $ constructed above, $ N=N(\ga)
$. Then $ R(x_j) \le \(\frac32\)^j R(x) \le 2^j R(x) $, and
\begin{multline*}
|x-y| \le \sum_{j=0}^{N-1} |x_j-x_{j+1}| \le
 \frac12 \sum_{j=0}^{N-1} R(x_j) \\
\le \frac12 R(x) \sum_{j=0}^{N-1} 2^j < 2^{N-1} R(x)
 \stackrel{\ref{1.3}}< 2^{3(\rho(x,y)+\epsilon)} R(x)\,.
\end{multline*}
\end{proof}

\section{Whitney-type partitions of unity}

For a smooth function $ f \colon \R^d \to \R $ we denote by $ |\nabla f(x)|
$ a norm of the gradient vector, say, $ |\nabla f(x)| = \sum_k \big|
\frac{\pd}{\pd x_k} f(x_1,\dots,x_d) \big| $ (the choice of the norm
does not matter), and by $ |\nabla^2 f(x)| $ a norm of the matrix of
second derivatives, say, $ |\nabla^2 f(x)| = \sum_{k,l} \big|
\frac{\pd^2}{\pd x_k \pd x_l} f(x_1,\dots,x_d) \big| $.

\begin{theorem}\label{2.1}
Let $ \rho $ be a special metric corresponding to $ R $, and $ U
\subset \R^d $ be a closed set. Then there exist a $ C^2 $-function $
f \colon \R^d \to [0,1] $ and $ \Const $ such that
\begin{gather*}
f(x) = 1 \quad \text{for all } x \in U \, , \\
f(x) = 0 \quad \text{for all } x \in \R^d \setminus U_{+4} \, , \\
\int \( R |\nabla f| + R^2 |\nabla^2 f| \) \, dm \le \Const \cdot \!
\int_{U_{+4}\setminus U} f \, dm \, .
\end{gather*}
\end{theorem}
\noindent Here and henceforth $ U_{+r} = \{ y \colon \exists x \in U \;
\rho(x,y) \le r \} $ is the $ r $-neighborhood of $ U $ with
respect to $ \rho $.

We denote for convenience $ B(x) = B \( x; \frac12 R(x) \) $ and $
\frac12 B(x) = B \( x; \frac14 R(x) \) $. By Lemma \ref{1.4},
\begin{equation}\label{2.15}
\forall x \;\; \forall y,z \in B(x) \quad \rho(y,z) \le 2 \, .
\end{equation}
The following fact follows immediately from \cite[Lemma
1.4.9]{Hormander}.

\begin{lemma}\label{2.2}
There exist a countable locally finite set $ S \subset \R^d $ and
$ \Const $ such that

\noindent\textup{(a)}
the balls $ \{ \frac12 B(s) \colon s \in S \} $ cover $ \R^d $,

\noindent\textup{(b)}
the multiplicity of the covering by twice larger balls
$ \{ B(s) \colon s \in S \} $ does not exceed $ \Const $.
\end{lemma}

The next lemma is essentially Theorem~1.4.10 from
\cite{Hormander}.
\begin{lemma}\label{2.3}
There exist constants $ \const $, $ \Const $, and $C^2$-functions
$ f_s \colon \R^d \to [0,1] $ for $ s \in S $ \textup{(}where $ S $ is
given by Lemma \textup{\ref{2.2})} such that

\noindent\textup{(a)}
$ f_s(x) = 0 $, unless $ x\in B(s) $;

\noindent\textup{(b)}
for all $ x \in \R^d $,
\[
\sum_{s\in S} f_s(x) = 1 \, ;
\]

\noindent\textup{(c)}
for all $ s \in S $,
\[
\int f_s \, dm \ge \const \cdot \, R^d(s) \, ;
\]

\noindent\textup{(d)}
for all $ s \in S $,
\[
\sup_{x\in\R^d} |\nabla f_s(x) | \le \frac{\Const}{R(s)} \, ,
\qquad
\sup_{x\in\R^d} |\nabla^2 f_s(x) | \le \frac{\Const}{R^2(s)} \, .
\]
\end{lemma}

\begin{proof}
We start with smooth functions $ g_s\colon \R^d \to [0,1] $ that satisfy
(a, d) and $ g_s(x) = 1 $ whenever $ x \in \frac12 B(s) $; the
latter implies (c) (for $ g_s $). Their sum
\[
g=\sum_{s\in S} g_s
\]
satisfies, for all $ x \in \R^d $,
\[
\const \le g(x) \le \Const \, .
\]
Indeed, the multiplicity of the covering (see Lemma \ref{2.2}) is an
upper bound; the lower bound (just $ 1 $) follows from \ref{2.2}(a).

It follows from condition (d) (for $g_s$) that
\[
| \nabla g(x) | \le \frac{\Const}{R(x)} \, , \qquad
| \nabla^2 g(x) | \le \frac{\Const}{R^2(x)}
\]
for all $ x $. It remains to take $ f_s = g_s / g $.
\end{proof}

\begin{proof}[Proof of Theorem \textup{\ref{2.1}}]
Lemmas \ref{2.2}, \ref{2.3} give us $ S $ and $ (f_s)_{s\in S} $; we
construct
\[
f = \sum_{s\colon B(s)\subset U_{+4}} f_s \, .
\]
By \ref{2.3}(a), $ f(x) = 0 $ for all $ x \notin U_{+4} $.
By \ref{2.3}(b), $ f(x) = 1 $ for all $ x \in U $ and moreover, for
all $ x \in U_{+2} $, since by \eqref{2.15}, $ B(s) \cap U_{+2} \ne
\emptyset $ implies $ B(s) \subset U_{+4} $.

We introduce a seminorm $ \|\cdot\| $,
\[
\| g \| = \int ( R |\nabla g| + R^2 |\nabla^2 g| ) \, dm
\]
for smooth functions $ g \colon \R^d \to \R $ such that this integral
converges. By \ref{2.3}(a,c,d),
\[
\| f_s \| \le \Const \cdot \! \int f_s \, dm
\]
for all $ s \in S $. Clearly, $ \| f \| \le \sum_{s\colon B(s)\subset
U_{+4}} \| f_s \| $, but moreover,
\[
\| f \| \le \sum_{s\colon B(s)\subset U_{+4} \setminus U} \| f_s \| \, ,
\]
since, taking into account that $ \nabla f = 0 $ outside $ U_{+4}
\setminus U_{+2} $, we have
\begin{multline*}
\| f \| = \int_{U_{+4}\setminus U_{+2}} ( R |\nabla f| + R^2 |\nabla^2
 f| ) \, dm \\
\le \sum_{s\colon B(s)\subset U_{+4}} \int_{U_{+4}\setminus
 U_{+2}} ( R |\nabla f_s| + R^2 |\nabla^2 f_s| ) \, dm \, ,
\end{multline*}
and the last integral vanishes for $ s $ such that $ B(s) \cap U \ne
\emptyset $. Finally,
\[
\| f \| \le \sum_{s\colon B(s)\subset U_{+4} \setminus U} \| f_s \| \le
\Const \cdot \sum_{s\colon B(s)\subset U_{+4} \setminus U} \int f_s \, dm
\le \Const \cdot \int_{U_{+4} \setminus U} f \, dm \, .
\]
\end{proof}

\section{The main lemma}

\begin{mainlemma}\label{3.1} There exists $ \const $ such that if a $ C^2
$-function $ u\colon \R^d \to \R $ and a special metric $ \rho $
\textup{(}corresponding to $ R $\textup{)} satisfy
\begin{gather}
\label{3.a} |u(x)| \le \const \cdot \, R^2 (x) \, , \\
\label{3.b} \De u (x) \ge -1
\end{gather}
for all $ x $, then
\[
\int_U \De u \, dm \le m ( U_{+4} \setminus U ) \, , \qquad -
\int_{U_{+4}} \De u \, dm \le m ( U_{+4} \setminus U )
\]
for every compact set $ U \subset \R^d $.
\end{mainlemma}

\noindent Still, $ U_{+4} = \{ y \colon \exists x \in U \; \rho(x,y)
\le 4 \} $.

\begin{proof}
Theorem \ref{2.1} gives us a function $ f \colon \R^d \to [0,1] $ that
equals $ 1 $ on $ U $, $ 0 $ outside $ U_{+4} $, and satisfies $
\int R^2 |\De f| \, dm \le \Const \cdot \int_{U_{+4}\setminus U} f
\, dm $. We have
\begin{align*}
-\int_{U_{+4}} \De u \, dm &= -\int f \De u \, dm -
 \int_{U_{+4}\setminus U} (1-f) \De u \, dm \, ; \\
-\int f \De u \, dm &= -\int u \De f \, dm \le \int |u| |\De f|
 \, dm \\
 &\le \const \cdot \int R^2 |\De f| \, dm \le \underbrace{
 \const \cdot \Const}_{\le 1} \cdot \int_{U_{+4}\setminus U} f \,
 dm \, ; \\
\int_{U_{+4}\setminus U} (1-f) (-\De u) \, & dm \le
 \int_{U_{+4}\setminus U} (1-f) \, dm \, ;
\end{align*}
thus,
\[
-\int_{U_{+4}} \De u \, dm \le \int_{U_{+4}\setminus U} f \, dm +
\int_{U_{+4}\setminus U} (1-f) \, dm = m ( U_{+4} \setminus U ) \,
.
\]

In order to prove the other inequality, we apply Theorem \ref{2.1}
to the closed set $ \R^d \setminus U_{+4} $ in place of $ U $ and
take $ 1-f $ in place of $ f $. This gives us $ f \colon \R^d \to [0,1]
$ that equals $ 1 $ on $ U $, $ 0 $ outside $ U_{+4} $, and
satisfies $ \int R^2 |\De f| \, dm \le \Const \cdot
\int_{U_{+4}\setminus U} (1-f) \, dm $. We have
\begin{align*}
\int_U \De u \, dm &= \int f \De u\, dm - \int_{U_{+4}\setminus
 U} f \De u \, dm \, ; \\
\int f \De u \, dm &= \int u \De f \, dm \le \int |u| |\De f|
 \, dm \\
&\le \const \cdot \int R^2 |\De f| \, dm \le
 \int_{U_{+4}\setminus U} (1-f) \, dm \, ; \\
\int_{U_{+4}\setminus U} f \cdot (-\De u) \, dm &\le
 \int_{U_{+4}\setminus U} f \, dm \, ;
\end{align*}
thus,
\[
\int_U \De u \, dm \le \int_{U_{+4}\setminus U} (1-f) \, dm +
\int_{U_{+4}\setminus U} f \, dm = m ( U_{+4} \setminus U ) \, .
\]
\end{proof}

\section{Tying zeroes of an entire function to the lattice points}

We fix once and forever a compactly supported smooth function $
\chi \colon \C \to [0,\infty) $ such that $ \int \chi \, dm = 1 $ and $
\chi(-z) = \chi(z) $ for all $ z $. The choice of $ \chi $ influences
only constants. Given an entire function $ \psi \colon \C \to \C $ and a
constant ($ \Const $), we define a function $ R \colon \C \to [1,\infty] $ by
\begin{equation}\label{4.4}
R(z) = \max_w \( \sqrt{ \Const \cdot \, ( 1 + (|\phi|*\chi)(w) ) }
- |w-z| \) \, ;
\end{equation}
here, as before,
\begin{equation}\label{4.5}
\phi(z) = \tfrac12 \log |\psi(z)| - \tfrac14 |z|^2 \, ,
\end{equation}
$ |\phi| $ is the pointwise absolute value of $ \phi $, and $
|\phi| * \chi $ is the convolution of these two functions.

Clearly, $ R(z) \ge \Const^{1/2} > 0 $ for all $ z $ (this time,
`$ \Const $' is the same as in \eqref{4.4}), and the function $ R
$ is a $ \Lip(1) $-function since it is an upper envelope of $
\Lip(1) $-functions. Therefore we may construct a special metric $
\rho $ corresponding to $ R $. Note that $ \rho(x,y) \le
\Const^{-1/2} |x-y| $ for all $ x,y $. If we replace `$ \Const $'
with `$ 4\Const $', we get another function $ R_2 $ such that $
R_2(z) \ge 2R(z) $ for all $ z $, and another special metric $
\rho_2 $ such that $ \rho_2(x,y) \le \frac12 \rho(x,y) $ for all $
x,y $.

Due to the Lipschitz property, the function $ R $ is finite
everywhere in $ \C $ provided that $ R(0) < \infty $.

\begin{theorem}\label{4.1}
There exists $ \Const $ in \eqref{4.4} with the following
property. For every entire function $ \psi $ satisfying $ R(0) <
\infty $ there exists a bijection between the lattice $ \sqrt\pi
\Z^2 $ and the zero set $ \psi^{-1}(0) $ \textup{(}counting with
multiplicities\textup{)} such that for every pair of corresponding
points $ z \in \psi^{-1}(0) $, $ \sqrt\pi(k+li) \in \sqrt\pi \Z^2
$
\begin{equation}\label{4.3}
| z - \sqrt\pi(k+li) | \le R \( \sqrt\pi(k+li) \) \, .
\end{equation}
\end{theorem}

\begin{proof}
In order to get \eqref{4.3} we will show that
\begin{equation}\label{4.9}
\rho ( z, \sqrt\pi(k+li) ) \le 5 \, ,
\end{equation}
where $ \rho $ is the special metric corresponding to $ R $.
In combination with Lemma \ref{1.5} it implies $ | z -
\sqrt\pi(k+li) | \le 2^{15} R \( \sqrt\pi(k+li) \) $; the constant
$ 2^{15} $ need not appear in \eqref{4.3}, since it can be
absorbed by $ \Const $ in \eqref{4.4}. Of course, there is nothing
sacred in the constant `5' on the RHS of \eqref{4.9}; any other
constant does the job as well.

Existence of a bijection between the lattice points and the set of
zeroes satisfying \eqref{4.9} follows from inequalities (to be
proven)
\begin{equation}\label{4.10}
n_\psi (U) \le n(U_{+5}) \, , \qquad n(U) \le n_\psi(U_{+5})
\end{equation}
for every compact $ U \subset \C $; here $ U_{+5} = \{ y \colon \exists
x \in U \; \rho(x,y) \le 5 \} $, $ n_\psi $ is the counting
measure on the set of zeroes, and $ n $ is the counting measure on
the lattice. Indeed, we say that the marriage between the bride $
\sqrt{\pi} (k+li) $ and the fianc\'e $ z\in\psi^{-1}(0) $ is
possible if $\rho (z, \sqrt{\pi} (k+li) )\le 5 $. By the classical
marriage lemma (we use its extension due to M.~Hall
\cite[p.7]{Aharoni}), the first inequality in \eqref{4.10}
guarantees that each bride can find a fianc\'e, and the second
inequality in \eqref{4.10} yields that each fianc\'e can find a
bride. Then a version of the Cantor-F.~Bernstein theorem
\cite[Theorem~1.1]{Aharoni} gives a matching which covers every
bride and every fianc\'e.

We split \eqref{4.10} into
\begin{gather}
\pi n(U) \le m(U_{+1}) \, , \qquad m(U) \le \pi n(U_{+1}) \, ;
 \label{4.11} \\
\pi n_\psi(U) \le m(U_{+4}) \, , \qquad m(U) \le \pi n_\psi(U_{+4}) \,
; \label{4.12}
\end{gather}
here $ m $ is the Lebesgue measure on $ \C $. Clearly, \eqref{4.11}
and \eqref{4.12} together imply \eqref{4.10}.

Inequalities \eqref{4.11} are easy to check, taking into account
that by Lemma~\ref{1.4} the $ \rho $-neighborhood $ U_{+1} $ of $
U $ contains a sufficiently large Euclidean neighborhood of $ U $,
provided that $ \Const $ in \eqref{4.4} is large enough. It
remains to prove \eqref{4.12}.

We will prove a seemingly weaker (but ultimately equivalent)
statement:
\begin{equation}\label{4.13}
\pi n_\psi (V) \le m (W) \, , \qquad m (V) \le \pi n_\psi (W)
\end{equation}
whenever $ V,W $ are such that $ V \subset U \subset U_{+4}
\subset W $, $ \One_U * \chi = 1 $ on $ V $, $ \One_{U_{+4}} *
\chi = 0 $ outside $ W $. (Here $ \One_U $ stands for the
indicator function of $ U $.) That is sufficient: having
\eqref{4.13} we get \eqref{4.12} after replacing `$ \Const $' with
`$ 4\Const $' in \eqref{4.4}, as follows. Considering the
corresponding $ R_2, \rho_2 $ and denoting the $ r $-neighborhood
w.r.t.\ $ \rho_2 $ by $ U^{+r} $ we have $ U^{+r} \supset U_{+2r}
$ (since $ \rho_2 \le \frac12 \rho $). Given $ \ti U $, we may
apply \eqref{4.13} to $ V = \ti U $, $ U = V_{+2} $, $ W = V^{+4}
\supset V_{+8} $ provided that $ \Const $ is large enough. We get
$ \pi n_\psi (\ti U) \le m (\ti U^{+4}) $, $ m (\ti U) \le \pi
n_\psi (\ti U^{+4}) $, which is \eqref{4.12}. It remains to prove
\eqref{4.13}.

Main lemma \ref{3.1} can be  applied to $ u = \phi * \chi $, since
condition \eqref{3.a} follows from \eqref{4.4}, and \eqref{3.b}
is checked readily:
\[
\De (\phi*\chi) = (\De \phi)*\chi = (\pi n_\psi - m)*\chi = \pi
n_\psi*\chi - m \ge - m.
\]
We get
\begin{gather*}
\int ( \One_U * \chi ) \De\phi \, dm = \int_U \De(\phi*\chi) \, dm \le
m ( U_{+4} \setminus U ) \, , \\
-\int ( \One_{U_{+4}} * \chi ) \De\phi \, dm = -\int_{U_{+4}}
\De(\phi*\chi) \, dm \le m ( U_{+4} \setminus U) \, .
\end{gather*}
However, $ \int ( \One_U * \chi ) \De\phi \, dm = \int ( \One_U * \chi
) \, d(\pi n_\psi - m) = \pi \int ( \One_U * \chi ) \, dn_\psi - m(U)
$, therefore
\[
\pi \int ( \One_U * \chi ) \, dn_\psi \le m(U_{+4}) \, , \qquad
m(U) \le \pi \int ( \One_{U_{+4}} * \chi ) \, dn_\psi \, ,
\]
and we get \eqref{4.13}:
\begin{gather*}
\pi n_\psi (V) \le \pi \int ( \One_U * \chi ) \, dn_\psi \le m(U_{+4})
 \le m(W) \, , \\
m(V) \le m(U) \le \pi \int ( \One_{U_{+4}} * \chi ) \, dn_\psi \le \pi
n_\psi (W) \, .
\end{gather*}
The proof of Theorem \ref{4.1} is completed.
\end{proof}

\section{Probabilistic arguments}

In the following lemma, $ (\Om,\F,P) $ is a probability space.

\begin{lemma}\label{5.1}
Let a random process $ \eta \colon \R^d \times \Om \to [0,\infty) $ satisfy
\begin{equation}\label{5.new1}
\Ex \exp \( c \eta(x) \) \le C
\end{equation}
for some $ C,c \in (0,\infty) $ and all $ x \in \R^d $,
and let $ R \colon \R^d \times \Om \to [0,\infty) $ be a random process
defined by
\[
R(x) = \max_y \( \sqrt{ \Const \cdot \, ( 1 + (\eta*\chi)(y) ) } -
|y-x| \) \, .
\]
Then there are constants $ c_1 $ and $ C_1 $ such that for all $ x \in
\R^d $,
\[
\Ex \exp \( c_1 R^2 (x) \) \le C_1 \, .
\]
\end{lemma}

\noindent The constants $ c_1 $ and $ C_1 $ depend on $ c $, $ C $, $
\Const $, the dimension $ d $, and the function $ \chi $ (introduced
in Sect.~4).

\begin{proof}
We will prove that
\begin{equation}\label{5.2}
\Pr{ R(0) > \la } \le C_1 \exp ( -c_1 \la^2 )
\end{equation}
for all $ \la $ large enough. It evidently implies $ \Ex \exp \( c_1
R^2(0) \) \le C_1 $ (with different $ c_1 $ and $ C_1 $). For other $
x $, $ R(x) $ is treated similarly.

We have
\[
\Pr{ R(0) > \la } = \Pr{ \exists x \; (\eta*\chi)(x) > \const \cdot \,
(\la + |x|)^2 - 1 } \, .
\]
For large $ \la $ we may discard `$ -1 $' on the RHS (at the expense
of changing the constant);
\[
\Pr{ R(0) > \la } \le \Pr{ \exists x \; (\eta*\chi)(x) > \const \cdot
\, \la^2 + \const \cdot \, |x|^2 } \, .
\]
Since the function $ \eta $ is non-negative, all values of the
convolution $ \eta * \chi $ are bounded by its values on a lattice:
given $ \chi $ there are constants $ c_2 $ and $ C_2 $ such that
\[
(\eta*\chi)(x) \le C_2 \max \{ (\eta*\chi)(y) \colon y \in c_2 \Z^d, \,
|y-x| \le C_2 \} \, .
\]
It follows that
\begin{align*}
\Pr{ R(0) > \la } &\le \Pr{ \exists x \in c_2 \Z^d \;
 (\eta*\chi)(x) > c_3 \la^2 + c_3 |x|^2 } \\
&\le \sum_{x \in c_2 \Z^d} \Pr{ (\eta*\chi)(x) > c_3 \la^2 + c_3 |x|^2
 } \, .
\end{align*}
Let $ c $ be a constant from \eqref{5.new1}. Then
\[
\Pr{ (\eta*\chi)(x) > c_3 \la^2 + c_3 |x|^2 }
 \le \frac{ \Ex \exp \( c (\eta * \chi)(x) \) }{ \exp ( c_4 \la^2 +
 c_4 |x|^2 ) }
\]
with $ c_4 = c_3 c $. However,
\[
\exp ( c (\eta * \chi) ) \le \( \exp ( c \eta ) \) * \chi
\]
(by convexity of `$ \exp $'), therefore
\[
\Ex \exp ( c (\eta * \chi) ) \le \( \Ex \exp ( c \eta ) \)
* \chi \stackrel{\eqref{5.new1}}\le C \, .
\]
Thus,
\[
\Pr{ (\eta*\chi)(x) > c_3 \la^2 + c_3 |x|^2 }
\le C \, \exp ( - c_4 \la^2 ) \exp ( - c_4 |x|^2 ) \, ,
\]
and the sum over $ x \in c_2 \Z^d $ does not exceed
$ C_1 \exp ( - c_4 \la^2 ) $, which proves \eqref{5.2}.
\end{proof}

Formula \eqref{0.0} defines a Gaussian random process $ \psi \colon \C
\times \Om_1 \to \C $ over the probability space $ (\Om_1,P_1) = \(
\C, \cN_\C(0,1) \)^{ \{0,1,2,\dots\} } $ (the space of coefficients $
\zeta_k $). The random variable $ \psi(0) = \zeta_0 $ is distributed $
\cN_\C(0,1) $. A simple exercise in integration shows that
\[
\Pr{ \big| \log |\zeta_0| \big| \ge s } = \Pr{ |\zeta_0| \ge e^s } +
\Pr{ |\zeta_0|\le e^{-s} } \le e^{-e^{2s}} + e^{-2s} \le 2e^{-2s} \, ,
\]
therefore $ \Ex \exp \( \big| \log |\psi(0)| \big| \) < \infty $.
We introduce a process $ \phi \colon \C \times \Om_1 \to \R $ by
\eqref{4.5}. The distribution of the random potential $ \phi(z) $
does not depend on $ z $ (see \cite{SoTsi}), therefore
\[
\Ex \exp ( 2|\phi(z)| ) \le \Const \, .
\]
By Lemma \ref{5.1} applied to $ \eta(z) = |\phi(z)| $,
\[
\Ex \exp \( \const \cdot \, R^2(z) \) \le \Const \, ,
\]
where $ R(\cdot) $ is a random process defined by \eqref{4.4}. It
follows that $ R(z) = O \( \sqrt{ \log |z| } \) $ for $ |z| \to
\infty $, almost surely. (The first part of the Borel-Cantelli
lemma gives it for $ z $ on a lattice; the Lipschitz property of $
R $ extends it to all $ z $.) By Theorem \ref{4.1}, for almost
every $ \om \in \Om_1 $ there exists a two-dimensional array $ \(
\xi_{k,l} (\om) \)_{k,l\in\Z} $ of complex numbers such that
\[
| \xi_{k,l} (\om) | \le R \( \sqrt\pi (k+li) \)
\]
for all $ k,l $, and the set $ \{ \sqrt\pi (k+li) + \xi_{k,l} \colon
k,l \in \Z \} $ is equal to the set $ \psi^{-1} (0) $ of zeroes of
$ \psi $. However, the array need not be unique.

\section{Final technicalities}

In order to get \emph{measurable} functions $ \om \mapsto
\xi_{k,l} (\om) $ one can use one of several well-known results about
measurable selectors, such as \cite[5.2.1, 5.2.5, 5.2.6, 5.4.3, 5.5.8,
5.7.1, 5.12.1]{Sri}. However, striving to keep the presentation
reasonably elementary, we borrow from \cite{Sri} only a special case
of Corollary 5.2.4, formulated below.

Given a complete separable metric space $ X $, we equip the set $
\K(X) $ of all compact subsets $ K \subset X $ with the Hausdorff
metric
\[
\dist (K_1,K_2) = \inf \{ \de > 0 \colon K_1 \subset (K_2)_{+\de}, \, K_2
\subset (K_1)_{+\de} \}
\]
and the corresponding Borel \sif\ (generated by open subsets of
$ \K(X) $). See the item `Spaces of compact sets' in
\cite[Sect.~2.4]{Sri}, see also the item `Effros Borel space' in
\cite[Sect.~3.3]{Sri}. The topology of $ \K(X) $ (known as the
Vietoris topology) is generated by sets of the following two forms:
\begin{gather*}
\{ K \in \K(X) \colon K \subset U \} \, , \\
\{ K \in \K(X) \colon K \cap U \ne \emptyset \} \, ,
\end{gather*}
where $ U $ runs over open subsets of $ X $. Unlike \cite{Sri}, we
treat the empty set as a point of $ \K(X) $; it is an isolated point
(take $ U = \emptyset $ in the first form above). The reader can check
that the same topology on $ \K(X) $ is generated also by functions $
\K(X) \to [-\infty,\infty) $ of the form
\[
K \mapsto \max_K \phi \, ,
\]
where $ \phi $ runs over bounded continuous functions $ X \to \R $;
for $ K = \emptyset $ the maximum is $ -\infty $.

\begin{theorem}\label{6.1}
There exists a Borel map $ s \colon \K(X) \setminus \{ \emptyset \} \to X $
such that $ s(K) \in K $ for all nonempty $ K \in \K(X) $.
\end{theorem}

A proof is given in \cite[5.2.4]{Sri}, but here is a hint: given $
\eps $, choose a countable $ \eps $-net $ \{ x_1, x_2, \dots \} $ of $
X $ and construct a Borel map $ s_\eps \colon \K(X) \setminus \{ \emptyset
\} \to \{ x_1, x_2, \dots \} $ such that $ \dist \( s_\eps (K), K \)
\le \eps $ for all $ K $.

\begin{lemma}\label{6.2}
Let $ X $ be a complete separable metric space, $ Y $ a metric space,
and $ f \colon X \to Y $ a continuous map. Then the map $ \K(X) \times Y
\to \K(X) $ defined by
\[
(K,y) \mapsto K \cap f^{-1} \( \{y\} \)
\]
is a Borel map.
\end{lemma}

\begin{proof}
Let $ \phi \colon X \to \R $ be a bounded continuous function; it is
sufficient to prove that the maximum of $ \phi $ on $ K \cap f^{-1} \(
\{y\} \) $ is a Borel function of $ (K,y) $. We use penalization:
\[
\max_{ x \in K \cap f^{-1} (\{y\}) } \phi(x) = \lim_{n\to\infty}
\max_{x\in K} \( \phi(x) - n \min(1,\dist(f(x),y)) \) \, .
\]
For each $ n $ the expression is continuous in $ y $ uniformly in $ K
$ and continuous in $ K $ for every $ y $, therefore it is continuous
in $ (K,y) $. The limit of a pointwise convergent sequence of such
functions is a Borel function.
\end{proof}

We apply the lemma to the separable Banach space $ X $ of all
two-dimensional arrays $ (\xi_{k,l}) $ of complex numbers satisfying $
|\xi_{k,l}| = o ( |k+li| ) $, with the norm
\[
\sup_{k,l} \frac{ |\xi_{k,l}| }{ 1 + | k+li | }
\]
(many other spaces could be used as well), the metrizable space $ Y $
of all locally finite measures on $ \C $, equipped with the topology
of local weak convergence, and the continuous map $ f \colon X \to Y $ that
sends $ (\xi_{k,l}) $ into the sum of unit-mass atoms at points $
\sqrt\pi (k+li) + \xi_{k,l} $. Every function $ R \colon \sqrt\pi \Z^2 \to
[0,\infty) $ such that $ R(z) = o(|z|) $ for $ |z| \to \infty $ leads
to a compact set $ K_R \subset X $,
\[
(\xi_{k,l}) \in K_R \iff \forall k,l \;\; |\xi_{k,l}| \le R \(
\sqrt\pi (k+li) \) \, .
\]
The map $ R \mapsto K_R $ is continuous w.r.t.\ the norm $ \| R \| =
\sup \( R(z) / (1+|z|) \) $.

As was shown in Sect.~5, almost each $ \om \in \Om_1 $ leads to an
entire function $ \psi_\om $, a measure $ n_\om = n_{\psi_\om} $
(recall \eqref{0.3}), and a function $ R_\om $ satisfying (much more
than) $ R_\om (z) = o (|z|) $ and such that the compact set
\[
K_\om = K_{R_\om} \cap f^{-1} \( \{ n_\om \} \)
\]
is nonempty. Measurability in $ \om $ of $ f_\om $ implies that of $
n_\om $, $ R_\om $ and, by Lemma \ref{6.2}, of $ K_\om $. Combined
with Theorem \ref{6.1}, it gives us the following result.

\begin{lemma}
There exist random variables $ \xi_{k,l} \colon \Om_1 \to \C $ such that
\begin{equation}\label{5.4}
\Ex \exp \( \const \cdot \, |\xi_{k,l}|^2 \) \le \Const
\end{equation}
and for almost all $ \om $ the set $ \{ \sqrt\pi (k+li) + \xi_{k,l} \colon
k,l \in \Z \} $ is equal to the set $ \psi^{-1} (0) $ of zeroes of $
\psi $.
\end{lemma}

The joint distribution of random variables $ \xi_{k,l} $ is a
probability measure on the space $ \Om = \C^{\Z^2} $ satisfying Items
(a), (b) of the main theorem. However, the measure need not be
shift-invariant.

The set of all probability measures on $ \Om $ satisfying both
Item (a) of the main theorem and \eqref{5.4} is convex, weakly
compact, and 
invariant under the action of $ \Z^2 $ (by continuous operators of
shift). By the Markov-Kakutani theorem \cite[Sect.~5.10.6]{DF}, the
action has a fixed point $ P $ in the set. This completes the
proof of the main theorem.

\bigskip
\filbreak
{
\small
\begin{sc}
\parindent=0pt\baselineskip=12pt
\parbox{2.3in}{
Mikhail Sodin\\
School of Mathematics\\
Tel Aviv University\\
Tel Aviv 69978, Israel
\smallskip
\emailwww{sodin@tau.ac.il}
{}
}
\hfill
\parbox{2.3in}{
Boris Tsirelson\\
School of Mathematics\\
Tel Aviv University\\
Tel Aviv 69978, Israel
\smallskip
\emailwww{tsirel@tau.ac.il}
{www.tau.ac.il/\textasciitilde tsirel/}
}
\end{sc}
}
\filbreak


\begin{thebibliography}{AAAA}

\bibitem{Aharoni} {\sc R. Aharoni},
Infinite matching theory,
Discrete Math. \textbf{95} (1991), 5--22.

\bibitem{AKT} {\sc M. Ajtai, J. Koml\'os and G. Tusn\'ady},
On optimal matchings,
Combinatorica \textbf{4} (1984), 259--264.

\bibitem{BSZ00a} {\sc P. Bleher, B. Shiffman and S. Zelditch},
Poincar\'e-Lelong approach to universality and
scaling of correlations between zeros,
Comm. Math. Phys. \textbf{208} (2000), 771--785.

\bibitem{DF} {\sc N. Dunford and J. Schwartz},
Linear operators. Part I.  General theory.
Interscience, N.Y., 1958.

\bibitem{Ha96} {\sc J.H. Hannay},
Chaotic analytic zero points: exact statistics for those of a
random spin state,
J. Phys. A: Math. Gen. \textbf{29} (1996), L101--L105.

\bibitem{HHP} {\sc C. Hoffman, A. E. Holroyd and Y. Peres},
A stable marriage of Poisson and Lebesgue.
In preparation.

\bibitem{HP} {\sc A. E. Holroyd and Y. Peres},
Extra heads and invariant allocations,
Ann. Prob. \textbf{33} (2005), 31--52.
\texttt{arXiv:math.PR/0306402}.

\bibitem{Hormander} {\sc L. H\"ormander},
The Analysis of Linear Partial Differential Operators, Vol. I.
Distribution Theory and Fourier Analysis, Springer Verlag,
Berlin, 1983.

\bibitem{Moser} {\sc J. Moser},
On the volume elements on a manifold,
Trans. Amer. Math. Soc. \textbf{120} (1965), 286--294.

\bibitem{PV} {\sc Y. Peres and B. Virag},
Zeros of the i.i.d. Gaussian power series: a conformally
invariant determinantal process,
Acta Mathematica \textbf{194} (2005), 1--35.
\texttt{arXiv:math.PR/0310297}

\bibitem{SoTsi} {\sc M. Sodin and B. Tsirelson},
Random complex zeroes, I. Asymptotic normality,
Israel Journal of Mathematics \textbf{144} (2004), 125--149.
\texttt{arXiv:math.CV/0210090}

\bibitem{SoTsi3} {\sc M. Sodin and B. Tsirelson},
Random complex zeroes, III. Decay of the hole probability.
Israel Journal of Mathematics \textbf{147} (2005), 371--379.
\texttt{arXiv:math.CV/0312258}

\bibitem{Sri} {\sc S. M. Srivastava},
A course on Borel sets, Springer Verlag, Berlin, 1998.

\bibitem{T} {\sc M. Talagrand},
Matching theorems and empirical discrepancy computations using
majorizing measures,
J. Amer. Math. Soc. \textbf{7} (1994), 455--537.

\end{thebibliography}
\end{document}